# ON THE LARGEST COMPONENT OF A RANDOM GRAPH WITH A SUBPOWER-LAW DEGREE SEQUENCE IN A SUBCRITICAL PHASE

BY B. G. PITTEL

*Ohio State University*

A uniformly random graph on $n$ vertices with a fixed degree sequence, obeying a $\gamma$ subpower law, is studied. It is shown that, for $\gamma > 3$, in a subcritical phase with high probability the largest component size does not exceed $n^{1/\gamma + \varepsilon_n}$, $\varepsilon_n = O(\ln \ln n / \ln n)$, $1/\gamma$ being the best power for this random graph. This is similar to the best possible $n^{1/(\gamma - 1)}$ bound for a different model of the random graph, one with independent vertex degrees, conjectured by Durrett, and proved recently by Janson.

**1. Introduction.** In a recently published book ([5], Section 1.2), Durrett formulated the following conjecture.

Let $\mathbf{p} = \{p_j\}_{j \geq 1}$ be a probability distribution. Let $D_1, \ldots, D_n$ be i.i.d. random variables, each having the distribution $\mathbf{p}$. Consider a graph on the vertex set $[n]$, chosen uniformly at random among all graphs with the degree sequence $(D_1, \ldots, D_n)$. For such a set of graphs to be nonempty, it is necessary that $\max D_i < n$ and $\sum_i D_i$ is even. (The first condition holds with probability approaching 1 if $\mathrm{E}[D] < \infty$, and the second condition holds with probability approaching 1/2 if $\mathrm{E}[D^2] < \infty$, and g.c.d.$\{j : p_j > 0\}$ is odd.) Durrett states that at a vicinity of a generic vertex $v$ the random graph looks like a tree rooted at $v$, and the number of direct descendants of every descendant of $v$ has a distribution $\mathbf{q} = \{q_j\}_{j \geq 0}$,

$$(1.1) \qquad q_j = \frac{(j+1)p_{j+1}}{\sum_{k \geq 1} k p_k}, \qquad j \geq 0.$$

---

Received August 2007; revised November 2007.

[1]Supported in part by NSF Grant DMS-04-06024.

If so, under the condition

$$\nu := \sum_j j q_j < 1,$$

one should expect that the component containing $v$, and even the largest component, are likely to be small compared to $n$. Specifically, Durrett conjectured that for the power-law distribution,

$$p_j = C j^{-\gamma}, \qquad \gamma > 3,$$

the likely size of the largest component should be of order $n^{1/(\gamma-1)}$, exactly. In other words, the largest component has size of order of the maximum vertex degree. Janson [6] has recently proved Durrett's conjecture.

In this paper we consider a different model of the random graph, in which a degree sequence is fixed. There is given a tuple $\mathbf{d} = \mathbf{d}(n) = (d_1, \ldots, d_n)$ of positive integers $d_1, \ldots, d_n < n$ such that $d_1 + \cdots + d_n$ is even. We consider a sample space $\mathcal{G}_{n,\mathbf{d}}$ of all graphs on $[n]$ with the degree sequence $\mathbf{d}$. Introduce the *empirical* degree distribution

$$\mathbf{p} = \{p_1, \ldots, p_{n-1}\}, \qquad p_j := \frac{|\{i \in [n] : d_i = j\}|}{n}.$$

Let $\mathbf{q} = \mathbf{q}(\mathbf{p})$ be defined by (1.1). Assuming that $\mathbf{p}$ obeys a subpower law, that is,

$$(1.2) \qquad\qquad p_j \le c j^{-\gamma}, \qquad 1 \le j \le n-1,$$

with $\gamma > 3$, we show that $\mathcal{G}_{n,\mathbf{d}}$ is nonempty. We prove that, under the condition

$$(1.3) \qquad\qquad \sum_{j \ge 1} j q_j \le 1 - \varepsilon, \qquad \varepsilon > 0,$$

the largest component in the graph $G_{n,\mathbf{d}}$, chosen uniformly at random from $\mathcal{G}_{n,\mathbf{d}}$, has size $C_n = O_p(n^{1/\gamma} \ln n)$, that is, $C_n/(n^{1/\gamma} \ln n)$ is bounded in probability. Similarly to Janson's result for the independent degrees model, the power $1/\gamma$ is the best possible for the fixed-degree-sequence model, since among the degree sequences $\mathbf{d}$ in question there are those with $\max_{v \in [n]} d_v$ of the exact order $n^{1/\gamma}$.

That, under the condition equivalent to (1.3), $C_n/n \to 0$ in probability, had already been proved by Molloy and Reed [7, 8]. They also proved that, under their form of the condition

$$\sum_{j \ge 1} j q_j \ge 1 + \varepsilon,$$

with high probability the random graph $G_{n,\mathbf{d}}$ has a giant component of size $\Theta(n)$, even being able to establish, under additional conditions, the limit of that size scaled by $n$.

Following the footsteps of Molloy and Reed, our proof is based on analysis of an algorithm that determines the component containing a given vertex. We construct a collection of exponential supermartingales in order to prove, via the optional sampling theorem, that the random growth of that component follows closely a certain deterministic path. See [1] and [9], where a similar approach was used for analysis of the site (bootstrap) and the bond percolation on a random *regular* graph.

**2. Main result and proofs.** Let $d_1, \ldots, d_n$ be positive integers, such that $d_1 + \cdots + d_n$ is even. Let $\mathcal{G}_{n,\mathbf{d}}$ denote the sample space of all graphs on the vertex set $[n]$ that have the degree sequence $\mathbf{d} = (d_1, \ldots, d_n)$. Denote by $G_{n,\mathbf{d}}$ the random graph which is distributed uniformly on $\mathcal{G}_{n,\mathbf{d}}$.

In parallel, let $M\mathcal{G}_{n,\mathbf{d}}$ denote the sample space of all multigraphs with the degree sequence $\mathbf{d}$. Let us describe the random multigraph $MG_{n,\mathbf{d}}$ suggested first by Bollobás [3]. Consider the disjoint sets $S_1, \ldots, S_n$ of cardinalities $d_1, \ldots, d_n$; set $S_i$ representing vertex $i \in [n]$. (Some people prefer assigning $d_i$ "half-edges" to a vertex $i \in [n]$, instead of sets $S_i$, but the difference is purely linguistic.) We know that $S := \bigcup_i S_i$ has an even cardinality $2m := d_1 + \cdots + d_n$. Introduce the sample space $\mathcal{P}_{n,\mathbf{d}}$ of all $(2m-1)!! = 1 \cdot 3 \cdots (2m-1)$ *pairings* on $S$. Let $P_{n,\mathbf{d}}$ be the random pairing distributed uniformly on $\mathcal{P}_{n,\mathbf{d}}$. Define $MG_{n,\mathbf{d}}$ as follows: two vertices $i, j \in [n]$ are joined by an edge iff there are $s' \in S_i$, $s'' \in S_j$ such that $\{s', s''\}$ is one of the pairs in $P_{n,\mathbf{d}}$. Obviously $MG_{n,\mathbf{d}}$ may well have loops and multiple edges. And it is *not* uniformly distributed on $M\mathcal{G}_{n,\mathbf{d}}$. However, conditioned on the event $A_n := \{\text{no loops, no multiple edges}\}$, $MG_{n,\mathbf{d}}$ is a simple graph distributed uniformly on $\mathcal{G}_{n,\mathbf{d}}$, hence can be viewed as the random graph $G_{n,\mathbf{d}}$. (This connection is due to the observation that every $G \in \mathcal{G}_{n,\mathbf{d}}$ induces the same number, $d_1! \cdots d_n!$, of pairings in $\mathcal{P}_{n,\mathbf{d}}$.)

Suppose that $\mathbf{d} = \mathbf{d}(n)$ is such that

$$(2.1) \qquad \liminf \mathrm{P}(A_n) > 0.$$

Under (2.1), any asymptotically rare (sure) event for $MG_{n,\mathbf{d}}$ is an asymptotically rare (sure) event for $G_{n,\mathbf{d}}$. And we will see that the probability estimates for the events in $M\mathcal{G}_{n,\mathbf{d}}$ become quite manageable once translated into the language of the space $\mathcal{P}_{n,\mathbf{d}}$.

Introduce

$$(2.2) \qquad \nu = \nu(n) = \frac{\sum_{i \in [n]} d_i(d_i - 1)}{\sum_{i \in [n]} d_i}.$$

$\nu$ can be interpreted as the expected *outdegree* of a nonroot vertex in a tree rooted at a given vertex $v$, which, heuristically, is how $G_{n,\mathbf{d}}$ looks like in a vicinity of $v$. Let

$$p_j = p_j(n) := \frac{1}{n} |\{i \in [n] : d_i = j\}|, \qquad j \in [n-1].$$

Then (2.2) becomes

$$\nu = \frac{\sum_{j \in [n-1]} j(j-1)p_j}{\sum_{j \in [n-1]} j p_j},$$

which is the ratio of the first two factorial moments of the distribution $\{p_j\}$. We denote the first moment, the average vertex degree, by $d = d(n)$.

We assume that $\mathbf{d} = \mathbf{d}(n)$ is such that

$$(2.3) \qquad \limsup_{n \to \infty} \sum_{j \in [n-1]} j^2 p_j < \infty, \qquad \lim_{n \to \infty} n^{-1} \sum_{j \in [n-1]} j^4 p_j = 0.$$

In fact, we assume a stronger condition, namely that $\{p_j\}$ is a subpower-law distribution, that is,

$$(2.4) \qquad\qquad p_j \le c\, j^{-\gamma}, \qquad \gamma > 3.$$

In this case, since

$$|\{i \in [n] : d_i = j\}| = n p_j \le c\frac{n}{j^\gamma} < 1 \qquad \forall j > j_n := \lfloor A(\gamma, c) n^{1/\gamma} \rfloor,$$

we see that

$$|\{i \in [n] : d_i = j\}| = 0, \qquad p_j = 0 \ \forall j > j_n.$$

In other words, $\max_{i \in [n]} d_i$, the largest vertex degree, is $j_n$, at most. That the first condition in (2.3) is met under (2.4) is obvious; the second condition holds true, since (2.4) implies

$$n^{-1} \sum_{j \in [n-1]} j^4 p_j = o(n^{-1/3});$$

see (2.19).

LEMMA. *Under the condition (2.3),*

$$\liminf \mathrm{P}(A_n) \ge \exp(-\hat{\nu}/2 - \hat{\nu}^2/4) > 0, \qquad \hat{\nu} := \limsup \nu.$$

NOTE. Applied to the degree sequence $\mathbf{d} = (d, \dots, d)$, this lemma yields a well-known asymptotic formula for the number of all $d$-regular graphs, due to Bender and Canfield [2].

We prove Lemma in the Appendix.

THEOREM. *Let $C_n$ denote the size of the largest component (cluster) of the random graph $G_{n,\mathbf{d}}$. Under the condition (2.4), for $\lambda = \lambda(n) \to \infty$ however slowly,*

$$\lim_{n \to \infty} \mathrm{P}\{C_n \le \lambda n^{1/\gamma} \ln n\} = 1,$$

*provided that*

(2.5) $$\limsup \nu < 1;$$

*in short, $C_n = O_P(n^{1/\gamma} \ln n)$.*

PROOF. We will prove the bound by upper-bounding the likely size of a component that contains a generic vertex $v \in [n]$. In view of Lemma, it suffices to bound the size of the component of the random multigraph $MG_{n,\mathbf{d}}$ that contains vertex $v$.

Notice that a subset $V$ of $[n]$ is the vertex set of this component iff for every $u \in V$ there exist $w_1, \ldots, w_k \in [n]$ such that, for some $s_0 \in S_v, s_1 \in S_{w_1}, \ldots, s_k \in S_{w_k}, s_{k+1} \in S_u$, all the pairs $\{s_0, s_1\}, \ldots, \{s_k, s_{k+1}\}$ are in $P_{n,\mathbf{d}}$. So we may, and will, deal with the corresponding "component" in $P_{n,\mathbf{d}}$ itself. We determine this component algorithmically, by adding to a current cluster of pairs exactly one new pair $\{s', s''\} \in P_{n,\mathbf{d}}$, where a *point* $s'$ is not in the current cluster of pairs, but has the same vertex "host" as one of the points in those pairs. [We will call them the (currently) active points.] If the point $s''$, the partner of the point $s'$, is hosted by a fresh vertex, $u$, then $u$ joins the current vertex cluster, and the $d_u - 1$ still unexplored points hosted by $u$ become active. As a result, the number of active points changes by $(-1) + (d_u - 1) = d_u - 2$. If $s''$ happens to be hosted by a vertex from the current vertex cluster, then the vertex cluster remains the same, but the number of active points decreases by 2.

Importantly, instead of generating the uniformly random paring $P_{n,\mathbf{d}}$ in advance, we can generate it one pair at a time, as called for by the algorithmic process. Namely, given a total ordering of the points in $S$, as $s'$ we pick the first, say, active point and pair it with a point $s''$, chosen uniformly at random among all points not in the pairs.

Let $A(t), I_j(t)$ denote the total number of the currently active points and the number of the currently inactive (not in the current cluster, i.e.) vertices after $t$ steps of the algorithm. In particular,

(2.6) $$A(0) = d_v, \qquad I_j(0) = np_j - \delta_{j,d_v}, \qquad j \leq j_n.$$

Introduce

$$I(t) = \sum_{j \leq j_n} j I_j(t),$$

the total number of inactive *points*. From the discussion above,

$$A(t+1) + I(t+1) = A(t) + I(t) - 2,$$

so that, by (2.6),

(2.7) $$A(t) + I(t) = A(0) + I(0) - 2t = nd - 2t,$$

where

$$d = n^{-1} \sum_{i \leq n} d_i = \sum_{j \leq j_n} j p_j$$

is the average vertex degree. From (2.7), the process will terminate no later than by time $t \leq nd/2$.

Clearly $\{A(t), \{I_j(t)\}_{j \leq j_n}\}_{t \geq 0}$ is a Markov chain, and if

$$t < T = T_v := \min\{\tau > 0 : \min\{A(\tau), I(\tau)\} = 0\},$$

then

$$(2.8) \qquad \mathrm{P}[I_j(t+1) = I_j(t) - 1 \mid \mathcal{F}_t] = -\frac{j I_j(t)}{A(t) + I(t) - 1},$$

$$\mathrm{E}[I_j(t+1) \mid \mathcal{F}_t] = I_j(t) + (-1)\frac{j I_j(t)}{A(t) + I(t) - 1},$$

$$(2.9) \qquad \mathrm{E}[A(t+1) \mid \mathcal{F}_t] = A(t) + (-2)\frac{A(t) - 1}{A(t) + I(t) - 1}$$

$$+ \sum_{j \geq 1} \frac{j I_j(t)}{A(t) + I(t) - 1}(j-2);$$

here $\mathrm{P}[\cdot \mid \mathcal{F}_t]$, $\mathrm{E}[\cdot \mid \mathcal{F}_t]$ denote the probability and the expectation conditioned on $\{A(t), \{I_j(t)\}_{j \leq j_n}\}$. Since $T$ is a stopping time, it follows from (2.8) and (2.7) that, for each $j \leq j_n$,

$$X_j(t) := \begin{cases} \dfrac{I_j(t)}{\prod_{\tau=0}^{t-1}(1 - j/(nd - 2\tau - 1))}, & t \leq T, \\ X_j(T), & t > T, \end{cases}$$

is a martingale, with

$$(2.10) \qquad \mathrm{E}[X_j(t)] \equiv X_j(0) = I_j(0) = np_j - \delta_{j,d_v}.$$

We want to show that, for $t = O(n^\alpha)$, $\alpha \in (\gamma^{-1}, 1 - \gamma^{-1})$, $X_j(t)$ is relatively close to $X_j(0)$, with probability very close to 1. (Of course, we focus on $\alpha$ close to $\gamma^{-1}$, since we expect the process to terminate around a time close to $n^{\gamma^{-1}}$.) To this end, first let us prove that the sequence

$$Q_j(t) := \exp[n^{\beta_j} X_j(t)/n], \qquad t \leq n^\alpha,$$

is "almost" a (super)martingale, provided that

$$(2.11) \qquad \gamma^{-1} + \alpha < 1,$$

$$2(\beta_j - 1) = \min\left\{0, -\alpha + (\gamma - 1)\frac{\ln j}{\ln n}\right\}.$$

Let $t < T$. Observe that, for $j \leq j_n$,

$$\prod_{\tau=0}^{t} \left(1 - \frac{j}{nd - 2\tau - 1}\right) = \exp\left(-\frac{j}{2} \int_{nd-2t}^{nd} \frac{dx}{x} + O(j/n)\right)$$

(2.12)
$$= \left(1 - \frac{2t}{nd}\right)^{j/2} (1 + O(j/n))$$

$$= 1 + O(jt/n)$$

$$= 1 + O(n^{\gamma^{-1}+\alpha-1}) \to 1.$$

Consequently, using

$$0 \leq I_j(t) - I_j(t+1) \leq 1, \qquad jI_j(t) \leq jI_j(0) = jp_j n \leq dn,$$

we obtain

$$n^{-1}X_j(t+1) - n^{-1}X_j(t)$$

$$= n^{-1} \frac{I_j(t+1) - I_j(t)}{\prod_{\tau=0}^{t}(1 - j/(nd - 2\tau - 1))}$$

(2.13)
$$\quad - n^{-1} \frac{I_j(t)j/(nd - 2t - 1)}{\prod_{\tau=0}^{t}(1 - j/(nd - 2\tau - 1))}$$

$$= O(n^{-1}|I_j(t+1) - I_j(t)|) + O(jI_j(t)n^{-2})$$

$$= O(n^{-1}).$$

Therefore, as $\beta_j \leq 1$,

$$\frac{Q_j(t+1)}{Q_j(t)} = \exp[n^{\beta_j}(n^{-1}X_j(t+1) - n^{-1}X_j(t))]$$

$$= 1 + n^{\beta_j - 1}(X_j(t+1) - X_j(t))$$

$$\quad + O(n^{2(\beta_j - 1)}(X_j(t+1) - X_j(t))^2).$$

Since $X_j(t)$ is a martingale, we have

$$E[X_j(t+1) - X_j(t) \mid \mathcal{F}_t] = 0.$$

Further, by (2.13) and (2.8),

$$E[(X_j(t+1) - X_j(t))^2 \mid \mathcal{F}_t]$$

$$= O(E[(I_j(t+1) - I_j(t))^2 \mid \mathcal{F}_t]) + O((jp_j)^2)$$

$$= O(jp_j + (jp_j)^2) = O(jp_j).$$

Consequently, for $t < T$, and trivially for $t \geq T$,

$$\frac{1}{Q_j(t)}\mathrm{E}[Q_j(t+1)\,|\,\mathcal{F}_t] = 1 + O(n^{2(\beta_j-1)}jp_j)$$

$$= 1 + O(n^{2(\beta_j-1)}j^{-\gamma+1})$$

$$= 1 + O(n^{-\alpha}),$$

the third equality and the fourth equality following from (2.4) and the definition of $\beta_j$ in (2.11), respectively. Thus, there exists $\varepsilon_n > 0$, $\varepsilon_n = O(n^{-\alpha})$, such that

$$\mathrm{E}[(1+\varepsilon_n)^{-t-1}Q_j(t+1)\,|\,\mathcal{F}_t] \leq (1+\varepsilon_n)^{-t}Q_j(t), \qquad t \leq n^{\alpha}.$$

It makes

$$\hat{Q}_j(t) := (1+\varepsilon_n)^{-t}Q_j(t), \qquad t \leq n^{\alpha},$$

a supermartingale, that differs from $Q_j(t)$ by a factor bounded away from both zero and infinity.

Given $j \leq j_n$, and $z > 0$, introduce a stopping time $\mathcal{T}_j(z)$, the first $t \leq n^{\alpha} \wedge T$ such that

$$\left| \frac{I_j(t)}{n}\prod_{\tau=0}^{t-1}\left(1 - \frac{j}{nd-2\tau-1}\right)^{-1} - \frac{I_j(0)}{n} \right| > \frac{z}{n^{\beta_j}},$$

and set $\mathcal{T}_j(z) = \lceil n^{\alpha}+1 \rceil$ if no such $t$ exists. By (2.10), for $t \leq n^{\alpha} \wedge T$ and $t < \min_j \mathcal{T}_j(z)$, we have

$$(2.14) \qquad I_j(t) = (np_j - \delta_{j,d_v})\prod_{\tau=0}^{t-1}\left(1 - \frac{j}{nd-2\tau-1}\right) + O(n^{1-\beta_j}z).$$

Applying the Optional Sampling Theorem to the supermartingale

$$\frac{\hat{Q}_j(t)}{\hat{Q}_j(0)} = (1+\varepsilon_n)^{-t}\exp[n^{\beta_j}(X_j(t)/n - X_j(0)/n)], \qquad t \leq n^{\alpha},$$

the stopping time $\mathcal{T}_j(z)$ (Durrett [4], Section 4.7), and using Markov inequality, we have: *uniformly* for $z > 0$, and $j \leq j_n$,

$$\mathrm{P}\{\mathcal{T}_j(z) = \lceil n^{\alpha}+1 \rceil\} = O(e^{-z}).$$

Choosing $z = \chi \ln n$, $(\chi > \gamma^{-1})$, and introducing

$$B_n = \left\{ \min_{j \leq j_n}\mathcal{T}_j(z) = \lceil n^{\alpha}+1 \rceil \right\},$$

we obtain then

$$(2.15) \quad \mathrm{P}(B_n) \geq 1 - O(j_n e^{-z}) = 1 - O(n^{\gamma^{-1}}e^{-\chi \ln n}) = 1 - O(n^{-\chi+\gamma^{-1}}).$$

Notice that, on the likely event $B_n$, (2.14) holds for all $t \leq n^\alpha \wedge T$.

Armed with (2.14), we turn our attention to (2.9) for

$$\mathrm{E}[A(t+1) \mid \mathcal{F}_t], \qquad t \leq n^\alpha \wedge T \text{ and } t < \min_j T_j(z).$$

By (2.14), for the sum in (2.10) we can write

$$
\begin{aligned}
(2.16) \quad \sum_{j \leq j_n} j(j-2) I_j(t) = {} & n \sum_{j \leq j_n} j(j-2) p_j \prod_{\tau=0}^{t-1} \left( 1 - \frac{j}{nd - 2\tau - 1} \right) \\
& + O(d_v^2) + O\left( \ln n \sum_{j \leq j_n} j^2 n^{1-\beta_j} \right).
\end{aligned}
$$

Here $d_v^2 = O(n^{2\gamma^{-1}})$, and by the definition of $\beta_j$,

$$
\begin{aligned}
\sum_{j \leq j_n} j^2 n^{1-\beta_j} & \leq \sum_{j \leq n^{\alpha/(\gamma-1)}} j^2 n^{1-\beta_j} + \sum_{j \leq j_n} j^2 \\
& = n^{\alpha/2} \sum_{j \leq n^{\alpha/(\gamma-1)}} \frac{1}{j^{(\gamma-5)/2}} + O(j_n^3) \\
& = \begin{cases} O(n^{3\alpha/(\gamma-1)} \ln n) + O(n^{3\gamma^{-1}}), & \gamma \leq 7, \\ O(n^{\alpha/2}) + O(n^{3\gamma^{-1}}), & \gamma > 7. \end{cases}
\end{aligned}
$$

Or

$$(2.17) \qquad \sum_{j \leq j_n} j^2 n^{1-\beta_j} = O(n^{3\gamma^{-1}}),$$

for $\alpha$ close to $\gamma^{-1}$. Furthermore, using (2.12),

$$
\begin{aligned}
(2.18) \quad & n \sum_{j \leq j_n} j(j-2) p_j \prod_{\tau=0}^{t-1} \left( 1 - \frac{j}{nd - 2\tau - 1} \right) \\
& = n \sum_{j \leq j_n} j(j-2) p_j \left( 1 - \frac{2t}{nd} \right)^{j/2} + O\left( \sum_{j \leq j_n} j^3 p_j \right).
\end{aligned}
$$

Here, by

$$1 - mx \leq (1-x)^m \leq 1 - mx + \binom{m}{2} x^2, \qquad x \geq 0,$$

we have

$$
\begin{aligned}
1 - \frac{2t}{nd} \frac{j}{2} & \leq \left( 1 - \frac{2t}{nd} \right)^{j/2} \\
& \leq 1 - \frac{2t}{nd} \frac{j}{2} + O(n^{-2} j^2 t^2);
\end{aligned}
$$

for $j = 1$ we need the lower bound as $j(j-2)|_{j=1} < 0$. Therefore

$$n \sum_{j \le j_n} j(j-2) p_j \left(1 - \frac{2t}{nd}\right)^{j/2}$$

(2.19)
$$\le n \sum_{j \le j_n} j(j-2) p_j + O\left(t \sum_{j \le j_n} j^3 p_j + n^{-1} t^2 \sum_{j \le j_n} j^4 p_j\right)$$

$$= n \sum_{j \le j_n} j(j-2) p_j + O(n^{\alpha+\gamma^{-1}} + n^{-1+2(\alpha+\gamma^{-1})})$$

$$= n \sum_{j \le j_n} j(j-2) p_j + O(n^{\alpha+\gamma^{-1}}),$$

as $\alpha + \gamma^{-1} < 1$; see (2.11). [We have used the bounds

(2.20)
$$\sum_{j \le j_n} j^3 p_j = O(n^{\max\{0,(4-\gamma)\gamma^{-1}\}} \ln n),$$

$$\sum_{j \le j_n} j^4 p_j = O(n^{\max\{0,(5-\gamma)\gamma^{-1}\}} \ln n),$$

which easily follow from (2.4).]

Combining (2.16)–(2.19), we obtain: for $t \le n^\alpha \wedge T$, $t < \min_j \mathcal{T}_j(z)$,

(2.21)
$$\sum_{j \le j_n} j(j-2) I_j(t) \le n \sum_{j \le j_n} j(j-2) p_j + O(n^{3\gamma^{-1}} \ln n),$$

if $\alpha$ is close to $\gamma^{-1}$. Notice that

$$\sum_j j(j-2) p_j = \left(\sum_j j p_j\right) \left(\frac{\sum_j j(j-1) p_j}{\sum_j j p_j} - 1\right) = d(\nu - 1),$$

so that

$$\limsup \sum_j j(j-2) p_j < 0,$$

which is the Molloy–Reed condition for the subcritical phase. So, by (2.21) and the condition $\gamma > 3$, (2.9) implies that

(2.22)
$$\mathrm{E}[A(t+1) \mid \mathcal{F}_t] \le A(t) - a \qquad \left[t \le T \wedge n^\alpha \text{ and } t < \min_j \mathcal{T}_j(z)\right];$$

$$a := \tfrac{1}{2} \limsup(1 - \nu) > 0,$$

for all $n$ large enough.

The rest is short. Set

$$A(t+1) = A(t) - a \qquad \left[t > T \wedge n^\alpha \text{ or } t \ge \min_j \mathcal{T}_j(z)\right].$$

Clearly the extended sequence $\{A(t)\}$ satisfies (2.22) for all $t$. Besides, since $T = T_v$ is the first time $\tau$ when

$$\min\{A(\tau), I(\tau)\} = 0,$$

we have

(2.23)
$$\{n^\alpha < T\} \cap \left\{\min_j \mathcal{T}_j(z) = \lceil n^\alpha + 1\rceil\right\} = \{n^\alpha < T\} \cap B_n$$
$$\subseteq \{A(\lceil n^\alpha\rceil) > 0\}.$$

Furthermore, since the maximum vertex degree is $j_n$ at most,

$$A(0) = d_v \leq j_n, \qquad |A(t+1) - A(t)| \leq j_n; \qquad j_n = O(n^{\gamma^{-1}}).$$

Also, reading out the conditional distribution $P\{A(t+1) - A(t) = i \,|\, \mathcal{F}_t\}$ from (2.10), and using (2.20),

$$E[(A(t+1) - A(t))^2 \,|\, \mathcal{F}_t] \leq 4 + \frac{2}{d} \sum_{j \leq j_n} j^3 p_j = O(n^{\max\{0,(4-\gamma)\gamma^{-1}\}} \ln n),$$

if $t \leq T \wedge n^\alpha$ and $t < \min_j \mathcal{T}_j(z)$. And the bound holds trivially for the larger values of $t$. Then

$$E[\exp(n^{-\gamma^{-1}}(A(t+1) - A(t))) \,|\, \mathcal{F}_t]$$
$$= 1 + n^{-\gamma^{-1}} E[A(t+1) - A(t) \,|\, \mathcal{F}_t]$$
$$\quad + O(n^{-2\gamma^{-1}} E[(A(t+1) - A(t))^2 \,|\, \mathcal{F}_t])$$
$$\leq 1 - an^{-\gamma^{-1}} + O(n^{\max\{-2\gamma^{-1},(2-\gamma)\gamma^{-1}\}} \ln n)$$
$$\leq 1 - bn^{-\gamma^{-1}}, \qquad b < a,$$

since $\gamma > 3$. Therefore

$$E[\exp(n^{-\gamma^{-1}}(A(t+1) - A(t))) \,|\, \mathcal{F}_t] \leq \exp(-bn^{-\gamma^{-1}}),$$

and then

$$E[\exp(n^{-\gamma^{-1}} A(t))] = O(\exp(-tbn^{-\gamma^{-1}})).$$

Hence

$$P\{A(t) > 0\} = O(\exp(-tbn^{-\gamma^{-1}})).$$

In particular, choosing

(2.24)
$$\alpha = \gamma^{-1} + \frac{\ln \ln n}{\ln n} + \frac{\eta}{\ln n}, \qquad \eta > 0,$$

which certainly satisfies the inequality $\gamma^{-1} + \alpha < 1$ in (2.11) for $n \geq n(\eta)$, we obtain

$$\mathrm{P}\{A([n^\alpha]) > 0\} = O(n^{-b\eta}).$$

By (2.23), we have then

$$\mathrm{P}\{(n^\alpha < T_v) \cap B_n\} = O(n^{-b\eta}),$$

and combining this estimate with (2.15) we conclude: for any fixed $\chi > 0$ and $\eta > 0$,

$$\mathrm{P}\{n^\alpha < T_v\} = \mathrm{P}\{e^\eta n^{\gamma^{-1}} \ln n < T_v\} = O(n^{-\chi + \gamma^{-1}} + n^{-b\eta}).$$

Of course, a bounded constant factor implicit in the big-Oh notation depends on $\chi$ and $\eta$. Thus, given $K > 0$, there exists $L = L(K)$ such that

$$\mathrm{P}\{L n^{\gamma^{-1}} \ln n < T_v\} \leq n^{-K-1}, \qquad v \in [n],$$

whence

$$\mathrm{P}\left\{\max_{v \in [n]} T_v > L n^{\gamma^{-1}} \ln n\right\} \leq n^{-K}.$$

It remains to notice that the component containing the vertex $v$ has size $T_v$ at most, so that $C_n$, the size of the largest component, is $\max_v T_v$, at most.

NOTE.  If, instead of (2.24), we had set

$$\alpha = \gamma^{-1} + \frac{\omega(n)}{\ln n}, \qquad \omega(n) \to \infty, \qquad \omega(n) = o(\ln n),$$

we would have proved that

(2.25)           $$\mathrm{P}\{\omega(n) n^{\gamma^{-1}} < T_v\} = O(n^{-\chi + \gamma^{-1}} + e^{-b\omega(n)}),$$

so that $T_v = O_p(n^{\gamma^{-1}})$. However, the $e^{-b\omega(n)}$ term in (2.25) would not have allowed us to deduce that $\max_{v \in [n]} T_v = O_p(n^{\gamma^{-1}})$ as well.  □

## APPENDIX

PROOF OF LEMMA.  Our argument is patterned after Bollobás's proof [3] of a similar, but more general, result for $\max_i d_i = O(1)$.

Let $X_n$ and $Y_n$ denote the total number of loops and the total number of pairs of parallel pairs in the random pairing $P_{n,\mathbf{d}}$. We want to show that, for every fixed $k$ and $\ell$,

$$\mathrm{E}[(X_n)_k (Y_n)_\ell] \sim \left(\frac{\nu}{2}\right)^{k+2\ell}, \qquad n \to \infty,$$

where $(a)_b$ stands for the falling factorial $a(a-1)\cdots(a-b+1)$. This would imply that $X_n$ and $Y_n$ are asymptotically independent, and Poisson distributed, with parameter $\nu/2$ and $(\nu/2)^2$, respectively, and the statement would follow, since

$$\mathrm{P}(A_n) = \mathrm{P}\{X_n = 0, Y_n = 0\}.$$

Combinatorially, $(X_n)_k (Y_n)_\ell$ is the total number of samples, with order and without replacement, of $k$ loops and of $\ell$ pairs of parallel edges from the random pairing $P_{n,\mathbf{d}}$. Given any such sample, let $S_{i_1}, \ldots, S_{i_{k+2\ell}}$ be the ordered sequence of sets such that $S_{i_j}$, $j \leq k$, contains the $j$th loop, and, for $1 \leq t \leq \ell$, the $t$th pair of parallel pairs is $(s_1, s_2)$, $(s_3, s_4)$, where $s_1, s_2 \in S_{i_{k+2t-1}}$, $s_3, s_4 \in S_{i_{k+2t}}$. We write

$$\mathrm{E}[(X_n)_k (Y_n)_\ell] = E_1 + E_2.$$

Here $E_1$ is the expected number of the samples such that $i_1 \neq \cdots \neq i_{k+2\ell}$, and $E_2$ is the expected number of all other samples, when at least two indices among $i_j$, $1 \leq j \leq k+2\ell$, coincide. Then

$$E_1 = \frac{(nd-2k-4\ell-1)!!}{(nd-1)!!} \sum_{i_1 \neq \cdots \neq i_{k+2\ell}} \prod_{s=1}^{k+2\ell} \binom{d_{i_s}}{2}.$$

EXPLANATION. Let a sequence $i_1 \neq i_2 \neq \cdots \neq i_{k+2\ell}$ be given. From each set $S_{i_j}$ we choose two points, in $\prod_{s=1}^{k+2\ell} \binom{d_{i_s}}{2}$ ways overall. We pair two points from each $S_{i_j}$, $j \leq k$, thus forming $k$ loops. For each $t \in [1, \ell]$, we match two chosen points in $S_{k+2t-1}$ with two chosen points in $S_{k+2t}$, in $2^\ell$ ways overall, and then divide by $2^\ell$ to account for irrelevance of the order in which every two sets, $S_{k+2t-1}$ and $S_{k+2t}$, appear in the sequence $S_{i_{k+1}}, \ldots, S_{i_{k+2\ell}}$.

Introduce

$$\Sigma_1 = \sum_{i_1, \ldots, i_{k+2\ell}} \prod_{s=1}^{k+2\ell} \binom{d_{i_s}}{2},$$

$$\Sigma_2 = \binom{k+2\ell}{2} \sum_{i_1 = i_2, i_3, \ldots, i_{k+2\ell}} \prod_{s=1}^{k+2\ell} \binom{d_{i_s}}{2};$$

so $\Sigma_1$ is a counterpart of $E_1$, with the indices $i_1, \ldots, i_{k+2\ell}$ allowed to coincide, and $\Sigma_2$ is an upper bound of the total sum of terms in $\Sigma_1$, but not in $E_1$. Clearly then

$$\frac{1}{(nd)^{k+2\ell}}(\Sigma_1 - \Sigma_2) \lesssim E_1 \lesssim \frac{1}{(nd)^{k+2\ell}}\Sigma_1.$$

Further

$$\Sigma_1 = \left(\sum_i \binom{d_i}{2}\right)^{k+2\ell} = (nd)^{k+2\ell}\left(\frac{1}{2d}\sum_j j(j-1)p_j\right)^{k+2\ell},$$

so that

$$\frac{\Sigma_1}{(nd)^{k+2\ell}} = \left(\frac{\nu}{2}\right)^{k+2\ell}.$$

Next

$$\Sigma_2 = \binom{k+2\ell}{2}\left(\sum_i \binom{d_i}{2}^2\right)\left(\sum_{i'} \binom{d_{i'}}{2}\right)^{k+2\ell-2},$$

so that

$$\frac{\Sigma_2}{(nd)^{k+2\ell}} = O\left(n^{-2}\sum_i \binom{d_i}{2}^2 \left(\frac{\nu}{2}\right)^{k+2\ell-2}\right) \to 0,$$

since, by (2.20),

$$n^{-2}\sum_i d_i^4 = n^{-1}\sum_{j \le j_n} j^4 p_j = O(n^{-1/3}).$$

Therefore

$$E_1 \sim \left(\frac{\nu}{2}\right)^{k+2\ell}, \qquad n \to \infty.$$

Finally

$$E_2 \le \frac{(nd - 2k - 4\ell - 1)!!}{(nd-1)!!} \cdot \Sigma_2 = O\left(n^{-2}\sum_i d_i^4\right) \to 0.$$

Therefore

$$\mathrm{E}[(X_n)_k(Y_n)_\ell] = E_1 + O(E_2) \sim \left(\frac{\nu}{2}\right)^{k+2\ell}, \qquad n \to \infty. \qquad \square$$

**Acknowledgment.** The critical comments by the Associate Editor were very helpful.

DEPARTMENT OF MATHEMATICS
OHIO STATE UNIVERSITY
COLUMBUS, OHIO 43210
USA
E-MAIL: bgp@math.ohio-state.edu